\input amstex
\documentstyle{amsppt}
\magnification=\magstep1
\NoBlackBoxes\NoRunningHeads
\topmatter
\title  
$C^*$-algebras arising from  
Dyck systems of topological Markov chains
\endtitle
\author Kengo Matsumoto
\endauthor
\affil 
Department of Mathematical Sciences, \\
Yokohama City  University,\\
22-2 Seto, Kanazawa-ku, Yokohama, 236-0027 Japan
\endaffil
\abstract
Let $A$
be an $N \times N $ irreducible matrix with entries in $\{0,1\}$.
We define the topological Markov Dyck  shift $D_A$ 
to be a nonsofic subshift consisting of the 
$2N$ brackets $(_1,\dots,(_N,)_1,\dots,)_N$ with both standard bracket rule 
and Markov chain rule coming from $A$.
The subshift is regarded as a subshift defined by the canonical generators 
$S_1^*,\dots, S_N^*, S_1,\dots, S_N $
of the Cuntz-Krieger algebra ${\Cal O}_A$.
We construct an irreducible 
$\lambda$-graph system 
${{\frak L}^{Ch(D_A)}}$
that presents the subshift $D_A$
so that  
we have an associated  simple purely infinite $C^*$-algebra
${\Cal O}_{{\frak L}^{Ch(D_A)}}$. 
We prove that ${\Cal O}_{{\frak L}^{Ch(D_A)}}$
is a universal unique $C^*$-algebra 
subject to some operator relations among 
$2N$ generating partial isometries.
Some examples are presented such that they are not stably isomorphic to any  
Cuntz-Krieger algebra.
\endabstract

\endtopmatter

\def\Zp{{ {\Bbb Z}_+ }}

\def\Ker{{{\operatorname{Ker}}}}

\def\id{{{\operatorname{id}}}}
\def\OLF{{{\Cal O}_{{\frak L}^{Ch(D_F)}}}}
\def\OLA{{{\Cal O}_{{\frak L}^{Ch(D_A)}}}}
\def\OLN{{{\Cal O}_{{\frak L}^{Ch(D_N)}}}}
\def\ALA{{{\Cal A}_{{\frak L}^{Ch(D_A)}}}}
\def\LCHDA{{{{\frak L}^{Ch(D_A)}}}}
\def\LCHDF{{{{\frak L}^{Ch(D_F)}}}}
\def\LCHDN{{{{\frak L}^{Ch(D_N)}}}}
\def\LCHLA{{{{\frak L}^{Ch(\Lambda_A)}}}}
\def\LWA{{{{\frak L}^{W(\Lambda_A)}}}}


 Keywords: $C^*$-algebra, Cuntz-Krieger algebra, 
subshift, $\lambda$-graph system,
Dyck shift,  K-theory,

Mathematics Subject Classification 2000:
Primary 46L80; Secondary 46L55, 37B10.

\bigskip

\heading 1. Introduction
\endheading

Let $\Sigma$ be a finite set with its discrete topology, 
that is called an alphabet.
Each element of $\Sigma$ is called a symbol.
Let $\Sigma^{\Bbb Z}$ 
be the infinite product space 
$\prod_{i=-{\infty}}^{\infty}\Sigma_{i},$ 
where 
$\Sigma_{i} = \Sigma$,
 endowed with the product topology.
 The transformation $\sigma$ on $\Sigma^{\Bbb Z}$ 
given by 
$\sigma({(x_i)}_{i \in \Bbb Z}) = (x_{i+1})_{i \in \Bbb Z}$ 
is called the full shift over $\Sigma$.
 Let $\Lambda$ be a closed subset of 
 $\Sigma^{\Bbb Z}$ such that $\sigma(\Lambda) = \Lambda$. 
 The topological dynamical system 
  $(\Lambda, \sigma\vert_{\Lambda})$
is called a subshift or a symbolic dynamical system.
 It is written as $\Lambda$ for brevity.
   There is a class of subshifts called sofic shifts,
that contains the topological Markov shifts.
Sofic shifts are presented by 
labeled graphs that are
 called  $\lambda$-graphs.

In [Ma], 
the  author  
has introduced a notion of 
$\lambda$-graph system 
as a generalization of 
$\lambda$-graph.
A $\lambda$-graph system 
$ {\frak L} = (V,E,\lambda,\iota)$
consists of a vertex set 
$V = V_0\cup V_1\cup V_2\cup\cdots$, an edge set 
$E = E_{0,1}\cup E_{1,2}\cup E_{2,3}\cup\cdots$, 
a labeling map
$\lambda: E \rightarrow \Sigma$
and a surjective map
$\iota_{l,l+1}: V_{l+1} \rightarrow V_l$ for each
$l\in \Zp$, 
where $\Zp$ denotes the set of all 
nonnegative integers.
An edge $e \in E_{l,l+1}$ has its source vertex $s(e)$ in $V_l,$
its terminal vertex $t(e)$ in $V_{l+1}$ 
and its label $\lambda(e)$ in $\Sigma$.

The theory of symbolic dynamical system has  a close relationship to formal language theory.
In the theory of formal language, 
there is a class of universal languages due to W. Dyck.
The symbolic dynamics generated by the languages are called the Dyck shifts 
$D_N$ (cf. [ChS], [Kr],[Kr2],[Kr3]). 
They are nonsofic subshifts.   
  Its alphabet consists of the $2N$ brackets: 
$(_1, \dots,(_N, )_1, \dots,)_N$.
 The forbidden words consist of words that do not obey the standard bracket rules. 
In [KM],  
a $\lambda$-graph system ${\frak L}^{Ch(D_N)}$ 
that presents the subshift $D_N$
has been introduced.
The $\lambda$-graph system is called the Cantor horizon  
$\lambda$-graph system 
for the Dyck shift 
$D_N$.
The K-groups for ${\frak L}^{Ch(D_N)}$,
that are invariant under topological conjugacy of the subshift $D_N$,
have been calculated ([KM]).

In [Ma2], a nuclear $C^*$-algebra ${\Cal O}_{\frak L}$
associated with a $\lambda$-graph system $\frak L$ has been introduced.
The class of the $C^*$-algebras contain the class of 
the Cuntz-Krieger algebras.
They are universal unique concrete $C^*$-algebras generated 
by finite families of partial isometries and  sequences of projections subject to certain operator relations  encoded by structure of the $\lambda$-graph systems.
Its K-groups 
$
K_i({\Cal O}_{\frak L}), i=0,1
$
are realized as the K-groups 
of the $\lambda$-graph system 
$\frak L$.
Therefore the results of [KM] imply that  
the $C^*$-algebras 
${\Cal O}_{\frak L^{Ch(D_N)}}$
for
$ N=2,3,\dots$
are unital simple purely infinite whose K-groups are
 $$
K_0({\Cal O}_{\frak L^{Ch(D_N)}}) \cong {\Bbb Z}/N{\Bbb Z} \oplus C(\frak K,\Bbb Z),
\qquad
K_1({\Cal O}_{\frak L^{Ch(D_N)}}) \cong 0
$$
where $C(\frak K,\Bbb Z)$ denotes the abelian group of all integer valued continuous functions on a Cantor discontinuum $\frak K$
([KM;Corolllary 3.17]).
Let $u_1,\dots,u_N$ be 
the canonical generating isometries 
of the Cuntz algebra ${\Cal O}_N$ 
that satisfy the relations:
$
\sum_{j=1}^{N} u_ju_j^* = 1, \,     u_i^*u_i = 1 
$
for 
$ i=1,\dots N. 
$
Then the bracket rule of the symbols $(_1, \dots,(_N,\,  )_1,\dots, )_N$
of the Dyck shift $D_N$ 
may be interpreted as the relations
$$
u_i^* u_i = 1  \qquad \text{ and } \qquad
u_i^* u_j = 0 \quad \text{ for } i \ne j
$$
of the partial isometries  
$u_1^*, \dots, u_N^*, u_1, \dots, u_N$
in the $C^*$-algebra ${\Cal O}_N$ (cf. (2.1)).

In the present paper, 
we consider a generalization of Dyck shifts $D_N$ 
by using the canonical generators of Cuntz-Krieger algebras 
${\Cal O}_A$
 for  
$N \times N$ matrices $A$ with entries in $\{ 0,1 \}$.
The generalized Dyck shift is denoted by $D_A$ 
and called the topological Markov Dyck  shift for $A$ (cf. [HIK], [KM2]).
Let $\alpha_1, \dots,\alpha_N, \, \beta_1, \dots, \beta_N$
be the alphabet of $D_A$.
They correspond to the brackets 
$(_1, \dots,(_N,\,  )_1,\dots, )_N$
respectively.
Let $t_1, \dots, t_N$ be 
the canonical generating partial isometries 
of the Cuntz-Krieger  algebra ${\Cal O}_A$ 
that satisfy the relations:
$
\sum_{j=1}^{N} t_jt_j^* = 1, \,     
t_i^*t_i = \sum_{j=1}^N A(i,j)t_jt_j^*
$
for 
$ i=1,\dots, N. 
$
Consider the correspondence
$
\varphi(\alpha_i) = t_i^*,
\varphi(\beta_i) = t_i, i=1,\dots,N.
$
Then a word $w$ of $\{\alpha_1, \dots, \alpha_N, \beta_1,\dots,\beta_N \}$
is defined to be admissible for the subshift $D_A$ 
precisely  if
the correspnoding element to $w$ through $\varphi$ in ${\Cal O}_A$
is not zero.
Hence we may recognize  $D_A$ to be  the subshift defined by the canonical generators of the Cuntz-Krieger algebra ${\Cal O}_A$.  
The subshifts  $D_A$ are not sofic in general and reduced to the Dyck shifts if all entries of $A$ are $1$.

We consider the Cantor horizon 
$\lambda$-graph system $\LCHDA$
for the topological Markov Dyck  shift $D_A$. 
The $\lambda$-graph system will be proved to be 
$\lambda$-irreducible with $\lambda$-condition (I) in the sense of [Ma5]
if the matrix is irreducible with condition (I) in the sense of Cuntz-Krieger [CK].
Hence the  associated $C^*$-algebra $\OLA$ is simple and purely infinite
because of [Ma5;Theorem 3.9].
We will show :
\proclaim{Theorem 1.1}
Let $A$ be  an $N \times N$ matrix  with entries in $\{ 0,1 \}$.
Suppose that $A$ is irreducible with condition (I).
 The $C^*$-algebra $\OLA$
associated with the $\lambda$-graph system
$\LCHDA$ is unital, separable,  nuclear, simple and purely infinite.
It is the unique $C^*$-algebra generated by $2N$ partial isometries 
$S_i, T_i, i=1,\dots,N$  
subject to the following operator relations:
$$
\align
\sum_{j=1}^{N} & ( S_jS_j^* +   T_jT_j^* )   =  1, \tag 1.1 \\
\sum_{j=1}^{N} & S_j^*S_j      = 1, \tag 1.2 \\
T_i^*T_i & =  \sum_{j=1}^N A(i,j) S_j^*S_j, \qquad i=1,2,\dots, N,\tag 1.3 \\
E_{\mu_1\cdots \mu_k} & = \sum_{j=1}^N A(j,\mu_1) S_jS_j^*
E_{\mu_1\cdots \mu_k}
S_jS_j^* + T_{\mu_1}E_{\mu_2\cdots \mu_k}T_{\mu_1}^*, \qquad k >1 \tag 1.4   
\endalign
$$
where 
$E_{\mu_1\cdots \mu_k}
= S_{\mu_1}^*\cdots S_{\mu_k}^*S_{\mu_k}\cdots S_{\mu_1}$,
$(\mu_1,\cdots,\mu_k)\in \Lambda_A^*$ 
 the set of admissible words of the topological Markov shift $\Lambda_A$ defined by the matrix $A$.
\endproclaim

If all entries of $A$ is $1$, 
the $\lambda$-graph system $\LCHDA$ becomes $\LCHDN$
so that the $C^*$-algebra $\OLA$ goes to the algebra  $\OLN$.
 If $A$ is  the Fibonacci matrix 
$
F =
\left[\smallmatrix
1 & 1 \\
1 & 0\\
\endsmallmatrix
\right]
$,
the $C^*$-algebra $\OLF$ is simple and purely infinite.
Its K-groups  are 
$
K_0(\OLF)  \cong  {\Bbb Z} \oplus C({\frak K},\Bbb Z)^{\infty},
\,
K_1(\OLF) \cong 0
$
where $C({\frak K},\Bbb Z)^{\infty}$ 
denotes the countable infinite direct sum 
of the group $C({\frak K},\Bbb Z)$ (cf.[Ma7]).
In general, the $C^*$-algebra ${\Cal O}_{\frak L}$ 
associated with a $\lambda$-graph system $\frak L$  has infinite family of generators. 
The $C^*$-algebras $\OLN$, $\OLF$ are finitely generated,
and
its $K_0$-groups however are not finitely generated.
Therefore  the algebras 
$\OLN$, $\OLF$ are not semiprojective
whereas Cuntz algebras and Cuntz-Krieger algebras are semiprojective
(cf. [Bla], [Bla2], [Ma3]).

\medskip
The  author would like to express  his sincerely  thanks to
Wolfgang Krieger whose suggestions and discussions   
made it possible to present this paper.

\heading 2. The topological Markov Dyck  shifts
\endheading
Throughout this paper $N$ is a fixed positive integer larger than $1$.

We consider the Dyck shift $D_N$ with alphabet 
$\Sigma = \Sigma^- \cup \Sigma^+$
where
$\Sigma^- = \{ \alpha_1,\cdots,\alpha_N \},
\Sigma^+ = \{ \beta_1,\cdots,\beta_N \}.
$
The symbols 
$ \alpha_i, \beta_i$
correspond to 
the brackets
$(_i,  )_i$
respectively, and  has the relations
$$
\alpha_i \beta_j
=
\cases
 {\bold 1} & \text{ if } i=j,\\
 0  & \text{ otherwise} \\
\endcases 
\tag 2.1
$$
for $ i,j = 1,\dots,N$  (cf. [Kr2],[Kr3]).
A word 
$\gamma_1\cdots\gamma_n $ 
of $\Sigma$
is defined to be admissible for $D_N$ precisely if
$
\prod_{m=1}^{n} \gamma_m \ne 0,
$
where
$
\prod_{m=1}^{n} \gamma_m 
$
means the products $\gamma_1 \cdots \gamma_n$ 
obtained by applying (2.1).

Let $A =[A(i,j)]_{i,j=1,\dots,N}$
be an $N\times N$ matrix with entries in $\{0,1\}$.
Throughout this paper,
 $A$ is assumed to have no zero rows or columns.
Consider the Cuntz-Krieger algebra ${\Cal O}_A$ for the matrix $A$
that is the universal $C^*$-algebra generated by 
$N$ partial isometries $t_1,\dots,t_N$ subject to the following relations:
$$
\sum_{j=1}^N t_j t_j^* = 1, 
\qquad
t_i^* t_i = \sum_{j=1}^N A(i,j) t_jt_j^* \quad \text{ for } i = 1,\dots,N
\tag 2.2
$$
([CK]).
Define a correspondence 
$\varphi_A :\Sigma \longrightarrow \{t_1^*,\dots, t_N^*, t_1,\dots, t_N \}$
by setting
$$ 
\varphi_A(\alpha_i) = t_i^*,\qquad 
\varphi_A(\beta_i) = t_i  \quad \text{ for } i=1,\dots,N.
$$
We denote by $\Sigma^*$ the set of all words 
$\gamma_1\cdots \gamma_n$ of elements of $\Sigma$.
Define the set
$$
{\frak F}_A = \{ \gamma_1\cdots \gamma_n \in \Sigma^* \mid
\varphi_A(\gamma_1)\cdots \varphi_A( \gamma_n) = 0 \}.
$$
Let $D_A$ be the subshift over $\Sigma$ whose forbidden words are 
${\frak F}_A.$
The subshift is called the topological Markov Dyck shift defined by $A$.
These kinds of  subshifts have first appeared in [HIK] in semigroup setting
and in [KM2] in more general setting without using $C^*$-algebras.
If all entries of $A$  are $1$, 
the partial isometries 
$ 
\varphi_A(\alpha_1),\dots, \varphi_A(\alpha_N),
\varphi_A(\beta_1),\dots, \varphi_A(\beta_N) 
$
satisfy the same relations as (2.1)
so that 
the subshift $D_A$ 
becomes the Dyck shift $D_N$
with $2N$ bracket.
We note the fact that 
 $\alpha_i \beta_j\in {\frak F}_A$  if $i\ne j$,
 and
 $\alpha_{i_n}\cdots \alpha_{i_1} \in {\frak F}_A$  
if and only if 
$\beta_{i_1}\cdots \beta_{i_n} \in {\frak F}_A$.
Consider the following  two subsystems of $D_A$
$$
\align
D_A^+ & = \{ {(\gamma_i)}_{i \in \Bbb Z} \in D_A \mid
\gamma_i \in \Sigma^+ \text{ for all } i \in \Bbb Z \},\\ 
D_A^- & = \{ {(\gamma_i)}_{i \in \Bbb Z} \in D_A \mid
\gamma_i \in \Sigma^- \text{ for all } i \in \Bbb Z \}.\\ 
\endalign
$$
The subshift  
$D_A^+$ is identified with the topological Markov shift 
$$
\Lambda_A = \{ {(x_i)}_{i \in \Bbb Z}\in \{ 1,\dots,N \}^{\Bbb Z}
 \mid A(x_i,x_{i+1}) = 1, i \in \Bbb Z \}
$$ 
defined by the matrix $A$ and similarly
$D_A^-$ is identified with 
the topological Markov shift $\Lambda_{A^t}$ 
defined by the transposed matrix $A^t$ of $A$.
Hence the subshift $D_A$ is recognized to contain both the topological Markov shifts 
$\Lambda_A$ and $\Lambda_{A^t}$.

\proclaim{Proposition 2.1}
If $A$ satisfies condition (I) in the sense of Cuntz-Krieger [CK],
the subshift $D_A$ is not sofic.
\endproclaim
\demo{Proof}
Put
$
X_{D_A^+} = \{ (\gamma_i)_{i \in \Bbb N} 
\mid (\gamma_i)_{i \in \Bbb Z}\in D_A^+\}
$
and
$
X_{\Lambda_A} = \{ (x_i)_{i \in \Bbb N} 
\mid (x_i)_{i \in \Bbb Z}\in \Lambda_A \}.
$
Since $A$ satisfies condition (I), 
we can for each $n=1,2,\dots $ find an element
$x(n)=(\beta_{n(i)})_{i\in \Bbb N} \in X_{D_A^+},$
where 
$(n(i))_{i \in \Bbb N} \in X_{\Lambda_A}$, 
so that $x(n) \ne x(k)$ for $n\ne k$.
Let
$\Gamma^-(x(n))$ be the predecessor set of $x(n)$ in $D_A$, that is, 
$$
\Gamma^-(x(n)) = \{ (\dots, y_{-2}, y_{-1}, y_0) \mid 
(\dots, y_{-2},y_{-1},y_0, \beta_{n(1)},\beta_{n(2)},\dots )\in D_A \}.
$$
Let for each $n=1,2,\dots $ and each $i \in \Bbb N$,
$\alpha_{n(i)} \in \Sigma^-$
so that
$\alpha_{n(i)} \beta_{n(i)} = {\bold 1}.$
Then
$(\cdots, \alpha_{n(2)}, \alpha_{n(1)}) \in \Gamma^-(x(k))$
if and only if
$k=n$.
Thus the predecessor sets
$\Gamma^-(x(n)), n= 1,2,\dots$
are mutually distinct,
so 
$D_A$ is not sofic.
\qed
\enddemo
Hence most irreducible matrix $A$ yield  non Markov subshifts $D_A$.

A $\lambda$-graph system $\frak L$ is said to present a subshift $\Lambda$
if the set of all admissible words of $\Lambda$ coincides with the set of all 
finite labeled sequences appearing in concatenating edges of $\frak L$. 
There are many $\lambda$-graph systems that present a given subshift.
Among them the canonical $\lambda$-graph system is a generalization 
of the left-Krieger cover graph for a sofic shift, 
and 
which strong shift equivalence class is 
invariant under topological conjugacy of subshifts ([Ma]).
The canonical $\lambda$-graph systems ${\frak L}^{C(D_N)}$ 
for the Dyck shifts  $D_N$ together with its K-groups 
have been  calculated in [Ma4].
One however sees that the $\lambda$-graph systems
${\frak L}^{C(D_N)}$
are not irreducible, so that the resulting $C^*$-algebras
 ${\Cal O}_{{\frak L}^{C(D_N)}}$
are not simple.
The Cantor horizon $\lambda$-graph system 
${\frak L}^{Ch(\Lambda)}$
for a Cantor horizon subshift $\Lambda$ 
is an irreducible component of the canonical $\lambda$-graph system 
${\frak L}^{C(\Lambda)}$.
The Cantor horizon $\lambda$-graph systems ${\frak L}^{Ch(D_N)}$ 
for the Dyck shifts  $D_N$ together with its K-groups 
have been  calculated in [KM] (cf. [Ma6]).
The associated  $C^*$-algebras
${\Cal O}_{{\frak L}^{Ch(D_N)}}$
are simple purely infinite  such that 
$$
K_0({\Cal O}_{{\frak L}^{Ch(D_N)}}) 
\cong {\Bbb Z}/N{\Bbb Z} \oplus C({\frak K},{\Bbb Z}), \qquad 
K_1({\Cal O}_{{\frak L}^{Ch(D_N)}}) \cong 0.
$$

In this paper we will study the Cantor horizon $\lambda$-graph systems 
$\LCHDA$
for the topological Markov Dyck  shifts $D_A$ 
and its associated $C^*$-algebras
${\Cal O}_{{\frak L}^{Ch(D_A)}}$
for $N\times N$ matrices $A$ with entries in $\{0,1\}$.

We denote by 
$B_l(D_A)$
and 
$B_l(\Lambda_A)$
 the set of admissible words of length 
$l$ of $D_A$
and that of 
$\Lambda_A$ respectively.
Let $m(l)$ be the cardinal number of $B_l(\Lambda_A)$.
We use lexcographic order from the left on the words of $B_l(\Lambda_A)$,
so that we may assign to a word $\mu_1\cdots \mu_l\in B_l(\Lambda_A)$ 
the number $N(\mu_1\cdots \mu_l)$ from $1$ to $m(l)$.
For example, 
if 
$A = 
\left[\smallmatrix 
1 & 1 \\
1 & 0 \\
\endsmallmatrix
\right],
$
then 
$$
\align
B_1(\Lambda_A)&  = \{ 1,2\},\qquad N(1) =1, \, N(2) = 2, \\
B_2(\Lambda_A)& = \{ 11, 12, 21 \},\qquad N(11) =1,\,  N(12) = 2,\,  N(21) = 3,
\endalign
$$
and so on.
Hence the set $B_l(\Lambda_A)$  bijectively corresponds 
to the set of natural numbers less than or equal to $m(l)$.  
Let us now describe 
the Cantor horizon $\lambda$-graph system $\LCHDA$ of $D_A$.
The vertices $V_l$ at level $l$ for $l \in \Zp$
are given by the admissible words of length $l$
consisting of the symbols of $\Sigma^+$.
We regard $V_0$ as a one point set of the empty word $\{ \emptyset \}$. 
Since $V_l$ is identified with $B_l(\Lambda_A)$,
we may write $V_l$ as
$$
V_l = \{ v^l_{N(\mu_1 \cdots \mu_l)}
 \mid \mu_1\cdots\mu_l\in B_l(\Lambda_A) \}.
$$
The mapping $\iota ( = \iota_{l,l+1}) :V_{l+1}\rightarrow V_l$ 
is defined by deleting the rightmost symbol of a corresponding word such as 
$$ 
\iota( v^{l+1}_{N(\mu_1\cdots \mu_{l+1})}) 
= v^l_{N(\mu_1 \cdots \mu_l)}\quad \text{ for }\quad
v^{l+1}_{N(\mu_1 \cdots \mu_{l+1})} \in V_{l+1}.
$$
We define an edge labeled $\alpha_j$ from
$v^l_{N(\mu_1 \cdots \mu_l)}\in V_l$ 
to
$v^{l+1}_{N(\mu_0\mu_1 \cdots \mu_l)}\in V_{l+1}$
precisely if
$\mu_0 = j,$
and 
 an edge labeled $\beta_j$ from
$v^l_{N(j\mu_1 \cdots \mu_{l-1})}
\in V_l$ 
to
$v^{l+1}_{N(\mu_1 \cdots \mu_{l+1})}\in V_{l+1}.$
For $l=0$, 
we define an edge labeled $\alpha_j$
form $v^0_1$ to $v^1_{N(j)}$, 
and
an edge labeled $\beta_j$
form $v^0_1$ to $v^1_{N(i)}$ if $A(j,i) =1$. 
We denote by $E_{l,l+1}$ the set of edges from $V_l$ to $V_{l+1}$.
Set $E = \cup_{l=0}^{\infty} E_{l,l+1}$.
It is easy to see that the resulting labeled Bratteli diagram with 
$\iota$-map
becomes a $\lambda$-graph system over $\Sigma$, that is denoted by
$\LCHDA$.

In the $\lambda$-graph system $\LCHDA$, 
we consider two $\lambda$-graph subsystems 
${\frak L}^{Ch(\Lambda_A)}$ and 
${\frak L}^{W(\Lambda_A)}$.
Both of the $\lambda$-graph subsystems have the same vertex sets as $\LCHDA$ together with the same $\iota$-maps as $\LCHDA$.
The edge set of $\LCHLA$ consist of edges labeled by $\Sigma^+$
in the edges of $\LCHDA$, 
whereas that of  $\LWA$ consist of edges labeled by $\Sigma^-$.
Hence $\LCHLA$ becomes a $\lambda$-graph system over $\Sigma^+$ and
 $\LWA$ becomes a $\lambda$-graph system over $\Sigma^-$. 
The latter $\lambda$-graph system 
is called the word $\lambda$-graph system in [KM2].
Since the union of the edge sets of $\LCHLA$ and $\LWA$ coincides with the edge set of $\LCHDA$, we may write $ \LCHDA$ as 
$$
       \LCHDA = \LCHLA \sqcup \LWA.
$$
We will prove that the $\lambda$-graph system 
$\LCHDA $ presents the subshift $D_A$.
We denote by $D_A^*$ and $\Lambda_A^*$ the set of admissible words of the subshifts $D_A$ and $\Lambda_A$ respectively.
\proclaim{Lemma 2.2}
For 
$\gamma_1\cdots \gamma_k \in D_A^*$ and $\mu_1\cdots \mu_l \in \Lambda_A^*$,
if the word 
$
\gamma_1\cdots \gamma_k \beta_{\mu_2}\cdots \beta_{\mu_l} 
$
is admissible in $D_A$, 
so is 
the word
$
\gamma_1\cdots \gamma_k \alpha_{\mu_1}\beta_{\mu_1}\beta_{\mu_2}\cdots \beta_{\mu_l}. 
$
\endproclaim
\demo{Proof}
As 
the word
$
\gamma_1\cdots \gamma_k \beta_{\mu_2}\cdots \beta_{\mu_l}
$
is admissible in $D_A$, one has
$$
\varphi_A(\gamma_1)\cdots \varphi_A(\gamma_k) 
t_{\mu_2}\cdots t_{\mu_l}t_{\mu_l}^*\cdots t_{\mu_2}^* \ne 0.
$$
By the condition
$\mu_1\cdots \mu_l \in \Lambda_A^*$ with
the relations (2.2), one sees
$$
t_{\mu_1}^*t_{\mu_1}t_{\mu_2}\cdots t_{\mu_l}t_{\mu_l}^*\cdots t_{\mu_2}^* 
=t_{\mu_2}\cdots t_{\mu_l}t_{\mu_l}^*\cdots t_{\mu_2}^*
$$
so that 
$$
\varphi_A(\gamma_1)\cdots \varphi_A(\gamma_k) t_{\mu_1}^*t_{\mu_1} 
t_{\mu_2}\cdots t_{\mu_l}t_{\mu_l}^*\cdots t_{\mu_2}^* \ne 0
$$
and hence the word
$
\gamma_1\cdots \gamma_k \alpha_{\mu_1}\beta_{\mu_1}\beta_{\mu_2}\cdots \beta_{\mu_l}
$ 
is admissible in $ D_A.$
\qed
\enddemo
For 
$\mu_1\cdots\mu_l\in B_l(\Lambda_A)$
and $k \le l$
we set 
$$
\Gamma^k_{D_A}(\beta_{\mu_1}\cdots\beta_{\mu_l})
 = \{ \gamma_1\cdots \gamma_k \in \Sigma^* \mid 
\gamma_1\cdots \gamma_k\beta_{\mu_1}\dots \beta_{\mu_l} \in D_A^* \}
$$
the $k$-predecessor set of the word $\beta_{\mu_1}\cdots\beta_{\mu_l}$ in $D_A$
and
$$
\align
& \Gamma^k_{\LCHDA}(v^l_{N(\mu_1\cdots \mu_l)}) \\
  = \{&  \gamma_1\cdots \gamma_k \in \Sigma^* \mid 
\text{there exist }e_i \in E , i=1,\dots,k \text{ such that }
 \gamma_i = \lambda(e_i) \\
 & \text{ for } i=1,\dots,k, \, 
 t(e_i) = s(e_{i+1}) \text{ for } i=1,\dots,k-1  \text{ and }
t(e_k) = v^l_{N(\mu_1\cdots \mu_l)} \}.
\endalign
$$
the $k$-predecessor set of the vertex $v^l_{N(\mu_1\cdots \mu_l)}$ in $\LCHDA$.
\proclaim{Lemma 2.3}
$
\Gamma^k_{D_A}(\beta_{\mu_1}\cdots\beta_{\mu_l})
=
\Gamma^k_{\LCHDA}(v^l_{N(\mu_1\cdots \mu_l)}).
$
\endproclaim
\demo{Proof}
We will prove the desired equality by induction on the length $k$.  

(1) Assume  that $k$ is $1$.

Take  $\mu_0 \in \{1,\dots,N\}$.
The symbol  
$\alpha_{\mu_0}$ belongs to 
$\Gamma^1_{D_A}(\beta_{\mu_1}\cdots\beta_{\mu_l})$
if and only if 
$\mu_0 = \mu_1$.
The latter condition is equivalent to the condition that 
$\alpha_{\mu_0}\in \Gamma^1_{\LCHDA}(v^l_{N(\mu_1\cdots \mu_l)}).$

The symbol  
$\beta_{\mu_0}$ belongs to
$\Gamma^1_{D_A}(\beta_{\mu_1}\cdots\beta_{\mu_l})$
if and only if 
$A(\mu_0,\mu_1) = 1$.
The latter condition is equivalent to the condition that 
$\beta_{\mu_0}\in \Gamma^1_{\LCHDA}(v^l_{N(\mu_1\cdots \mu_l)}).$

(2) Assume next that the desired equality holds for a fixed 
$k$ with $k+1 \le l$.
Take a word $\gamma_1\cdots\gamma_{k+1} \in \Sigma^*$.
We have two cases.

Case 1:  $\gamma_{k+1} = \alpha_{\mu_0}$ for some 
$\mu_0 \in \{1,\dots,N\}$.

Assume that
$
\gamma_1\cdots \gamma_{k+1} 
$ 
belongs to
$
\Gamma^{k+1}_{D_A}(\beta_{\mu_1}\cdots\beta_{\mu_l})
$
and hence 
$
\gamma_1 \cdots \gamma_k \alpha_{\mu_0}  \beta_{\mu_1} \cdots \beta_{\mu_l} 
$ 
is admissible in $ D_A$.
One then sees that 
$
\mu_0 = \mu_1$.
Since $t_{\mu_0}^* t_{\mu_0}$ is a projection in the algebra 
${\Cal O}_A$, the word      
$
 \gamma_1 \cdots \gamma_k \beta_{\mu_2} \cdots \beta_{\mu_l} 
$ 
is admissible in $ D_A$.
Hence
$$
\gamma_1\cdots \gamma_k \in
\Gamma^k_{D_A}(\beta_{\mu_2}\cdots\beta_{\mu_l}).
$$
By the hypothesis of induction, one has 
$$
\gamma_1\cdots \gamma_k \in
\Gamma^k_{\LCHDA}(v^{l-1}_{N(\mu_2\cdots \mu_l)}).
$$
Since $\mu_1 \mu_2\cdots \mu_l$ is admissible in $\Lambda_A$, 
there exists an edge $e \in E_{l-1,l}$ in $\LCHDA$ such that 
$\lambda(e) =  \alpha_{\mu_1}$ and 
$s(e) =v^{l-1}_{N(\mu_2\cdots \mu_l)},
t(e) =v^l_{N(\mu_1\cdots \mu_l)}$.
Hence we know that 
$$
\gamma_1\cdots \gamma_{k+1} \in
\Gamma^{k+1}_{\LCHDA}(v^l_{N(\mu_1\cdots \mu_l)}).
$$
Conversely assume that 
$
\gamma_1\cdots \gamma_{k+1}
$ belongs to
$
\Gamma^{k+1}_{\LCHDA}(v^l_{N(\mu_1\cdots \mu_l)}).
$
Since $\gamma_{k+1} = \alpha_{\mu_0}$,
one has 
that 
$
\mu_0 = \mu_1.
$
Hence
$$
\gamma_1\cdots \gamma_k \in
\Gamma^k_{\LCHDA}(v^{l-1}_{N(\mu_2\cdots \mu_l)}).
$$
By the hypothesis of induction, the word 
$
 \gamma_1 \cdots  \gamma_k \beta_{\mu_2} \cdots \beta_{\mu_l}
$ 
is admissible in $ D_A$.
Since $\mu_1 \mu_2\cdots \mu_l$ is admissible in $\Lambda_A$, 
by the preceding lemma,
$
 \gamma_1 \cdots  \gamma_k \alpha_{\mu_1} \beta_{\mu_1}  \beta_{\mu_2} \cdots \beta_{\mu_l} 
$ 
is admissible in $D_A$
so that 
$$
\gamma_1\cdots \gamma_{k+1}\in
\Gamma^{k+1}_{\LCHDA}(\beta_{\mu_1} \cdots \beta_{\mu_l}).
$$

Case 2: $\gamma_{k+1} = \beta_{\mu_0}$ for some 
$\mu_0 \in \{1,\dots,N\}$.

Assume that 
$
\gamma_1\cdots \gamma_{k+1}
$
belongs to
$
\Gamma^{k+1}_{D_A}(\beta_{\mu_1}\cdots\beta_{\mu_l}).
$
Then 
$$
\gamma_1\cdots \gamma_k\in
\Gamma^k_{D_A}(\beta_{\mu_0}\beta_{\mu_1}\cdots\beta_{\mu_{l-2}}).
$$
By the hypothesis of induction,
we have
$$
\gamma_1\cdots \gamma_k\in
\Gamma^k_{\LCHDA}(v^{l-1}_{N(\mu_0\cdots \mu_{l-2})}).
$$
Since $\mu_0 \mu_1\cdots \mu_l$ is admissible in $\Lambda_A$, 
there exists an edge $e \in E_{l-1,l}$ in $\LCHDA$ such that 
$\lambda(e) = \beta_{\mu_0}$ and 
$s(e) = v^{l-1}_{N(\mu_0\cdots \mu_{l-2})},
t(e) = v^l_{N(\mu_1\cdots \mu_l)}$.
Hence we know  
$$
\gamma_1\cdots \gamma_{k+1} \in
\Gamma^{k+1}_{\LCHDA}(v^l_{N(\mu_1\cdots \mu_l)}).
$$
Conversely assume that 
$
\gamma_1\cdots \gamma_{k+1}
$ belongs to
$
\Gamma^{k+1}_{\LCHDA}(v^l_{N(\mu_1\cdots \mu_l)}).
$
Hence
$$
\gamma_1\cdots \gamma_k \in
\Gamma^k_{\LCHDA}(v^{l-1}_{N(\mu_0\cdots \mu_{l-2})}).
$$
As
$
\Gamma^k_{\LCHDA}(v^{l-1}_{N(\mu_0\cdots \mu_{l-2})})
=
\Gamma^k_{\LCHDA}(v^{l+1}_{N(\mu_0\cdots \mu_l)}),
$
one has
by the hypothesis of induction
$$
\gamma_1\cdots \gamma_k \in
\Gamma^k_{D_A}(\beta_{\mu_0}\cdots \beta_{\mu_l}).
$$
Hence we have 
$$
\gamma_1\cdots \gamma_{k+1} \in
\Gamma^{k+1}_{D_A}(\beta_{\mu_1}\cdots\beta_{\mu_l}).
$$
Therefore the desired equality holds for all $k$ with $k\le l$.
\qed
\enddemo
\proclaim{Proposition 2.4}
 The $\lambda$-graph system $\LCHDA$ presents the subshift $D_A$.
 \endproclaim
\demo{Proof}
Put
$
X_{\Lambda_A} = \{ {(\mu_i )}_{i \in \Bbb N}  
\mid {(\mu_i)}_{i \in \Bbb Z}\in \Lambda_A\}
$
and
$
X_{D_A} = \{ {(\gamma_i)}_{i \in \Bbb N}  
\mid {(\gamma_i)}_{i \in \Bbb Z}\in D_A\}.
$
We faithfully represent the Cuntz-Krieger algebra ${\Cal O}_A$
on the Hilbert space $\frak H$ 
whose complete orthonormal basis are given by the vectors
$$
e_{\mu_1}\otimes e_{\mu_2}\otimes \cdots \qquad \text{ for }\quad (\mu_1, \mu_2,\dots  ) \in X_{\Lambda_A}
$$
by using the creation operators $t_i, \, i = 1,\dots,N$ 
on $\frak H$ defined by 
$$
t_i (e_{\mu_1} \otimes e_{\mu_2} \otimes \cdots )
=
\cases
e_i\otimes e_{\mu_1} \otimes e_{\mu_2} \otimes \cdots  & \text{ if } 
A(i,\mu_1) =1,\\
0 & \text{ otherwise. }
\endcases
$$
We may identify  $\varphi_A(\alpha_i)$ and $ \varphi_A(\beta_i) $
with  the operators 
$t_i^*$ and $t_i$ on $\frak H$ respectively. 
For a word $\gamma_1\cdots \gamma_k \in \Sigma^*$, 
it folllows that
$\gamma_1\cdots \gamma_k$ is admissible in $D_A$ 
if and only if 
there exists a sequence 
$(\mu_1,\mu_2, \dots ) \in X_{\Lambda_A}$
such that
$
\varphi_A(\gamma_1) \cdots \varphi_A(\gamma_k)
e_{\mu_1} \otimes e_{\mu_2} \otimes \cdots 
$
is a nonzero vector.
The latter condition is equivalent to the condition 
$
(\gamma_1,\dots, \gamma_k,\mu_1,\mu_2, \dots  ) \in X_{D_A}.
$
This is equivalent to the condition 
$
\gamma_1\cdots \gamma_k \in
\Gamma^k_{D_A}(\beta_{\mu_1}\beta_{\mu_2}\cdots\beta_{\mu_l})
$
for all $l \ge k$.
Therefore by the preceding lemma, 
the subshift 
$\Lambda_{\LCHDA}$ presented by the $\lambda$-graph system $\LCHDA$ is $D_A$.
\qed
\enddemo

Hence we automatically know that the $\lambda$-graph systems
$\LCHLA$ and $\LWA$ present the subshifts $D_A^+$ and $D_A^-$ respectively. 

A $\lambda$-graph system $\frak L$ satisfies $\lambda$-condition (I) if
for every vertex $v \in V_l$ of $\frak L$ there exist at least 
two paths with distinct label sequences starting with the vertex $v$ and terminating with the same vertex.

${\frak L}$ is said to be $\lambda$-irreducible if
for an ordered pair of vertices $u, v \in V_l,$ 
there exists a number $L_l(u,v) \in \Bbb N$ such that 
for a vertex $w \in V_{l+L_l(u,v)}$ 
with
$\iota^{L_l(u,v)}(w) = u,$
there exists a path $\xi$ in $\frak L$ 
such that 
$
s(\xi) = v, \, 
t(\xi) = w,
$
where $\iota^{L_l(u,v)}$ means the $L_l(u,v)$-times compositions of $\iota$, 
and $s(\xi), t(\xi)$ denote 
the source vertex, the terminal vertex of $\xi$ respectively ([Ma5]).
\proclaim{Proposition 2.5}
Let $A$ be an $N \times N$ matrix with entries in $\{0,1\}$.
\roster
\item"(i)" If $A$ satisfies condition (I) in the sense of Cuntz-Krieger [CK],
the $\lambda$-graph system $\LCHLA$ satisfies $\lambda$-condition (I).
\item"(ii)" If $A$ is irreducible,
the $\lambda$-graph system $\LCHLA$ is $\lambda$-irreducible.
\endroster
Hence if $A$ is an irreducible matrix with condition (I), then
both the $\lambda$-graph systems 
$\LCHLA$ and $\LCHDA$ are $\lambda$-irreducible with $\lambda$-condition (I).
\endproclaim
\demo{Proof}
(i)
Suppose that $A$ satisfies condition (I).
In the $\lambda$-graph system
${\frak L}^{Ch(\Lambda_A)}$, 
let $v_i^l$ be a vertex in $V_l$.
We write 
$i = N(i_1 \cdots i_l)$ 
for 
${i_1}\cdots {i_l} \in B_l(\Lambda_A)$.
By condition (I) for $A$,
there exist
$
\mu = \mu_1\cdots\mu_r, \,  
 \nu = \nu_1\cdots\nu_r \in B_r(\Lambda_A)$
 such that
 $\mu \ne \nu,  \mu_1 = \nu_1 = {i_l}$
 and 
 $\mu_r = \nu_r.$ 
Take
$
\eta_{r+1}\cdots\eta_{2l+2r-1}
 \in B_{2l+r-1}(\Lambda_A)$
such that
$
\mu_r\eta_{r+1}\dots\eta_{2l+2r-1}
 \in B_{2l+r}(\Lambda_A).
 $
We put
$
\mu_n =\nu_n= \eta_n
$
for
$ n= r+1,\dots, 2l+2r-1.
$
Put $L' = 2l+2r-1$.
Let $v_j^{L'} \in V_{L'}$ be the vertex in ${\frak L}^{Ch(\Lambda_A)}$
such that
$j =N(\mu_r \mu_{r+1}\cdots \mu_{2l+2r-2}) 
(=N(\nu_r \nu_{r+1}\cdots \nu_{2l+2r-2}))$.
Then there exist two paths labeled 
$\beta_{i_1}\cdots \beta_{i_l} \beta_{\mu_1}\cdots \beta_{\mu_{r-1}}$
and  
$\beta_{i_1}\cdots \beta_{i_l} \beta_{\nu_1}\cdots \beta_{\nu_{r-1}}$
whose sources are both $v_i^l$ and terminals are both $v_j^{L'}$.
Hence $\LCHLA$ satisfies $\lambda$-condition (I).

(ii)
In the $\lambda$-graph system
${\frak L}^{Ch(\Lambda_A)}$, 
let $v_i^l, v_j^l$ 
be  vertices in $V_l$.
We write 
$i=N(i_1\cdots i_l), j=N(j_1\cdots j_l)
$
for
$
i_1\cdots i_l, 
j_1\cdots j_l
 \in B_l(\Lambda_A)
 $
respectively.
As $\Lambda_A$ is irreducible,
there exits a word
$\eta_1 \cdots \eta_L \in B_L(\Lambda_A)$
such that
$
j_1\cdots j_l \eta_1 \cdots \eta_L i_1\cdots i_l \in B_{2l+L}(\Lambda_A).
$
We may assume 
$L \ge l$.
For $v_h^{2l+L} \in V_{2l+L}$ with
$\iota^{l+L}(v_h^{2l+L}) = v_i^l, h= 1,\dots,m(2l+L)$
we have
$
h = N(i_1\cdots i_l 
\mu_{l+1}\cdots \mu_{2l+L})
$
for some
$
\mu_{l+1} \cdots \mu_{2l+L}
\in B_{l+L}(\Lambda_A).
$
Then there exists a path labeled 
$ \beta_{j_1}\cdots \beta_{j_l} \beta_{\eta_1} 
\cdots \beta_{\eta_L}$
whose source is $v_j^l$ and whose terminal is $v_h^{2l+L}$.
This means 
that ${\frak L}^{Ch(\Lambda_A)}$
is $\lambda$-irreducible.
\qed
\enddemo

\heading 3. The $C^*$-algebra $\OLA$
\endheading
This section is devoted to studying operator relations among generators of 
 the algebra $\OLA$ to prove Theorem 1.1.
A general structure for the $C^*$-algebra ${\Cal O}_{\frak L}$
associated with a  $\lambda$-graph system $\frak L$
has been studied in [Ma2] as in the following way:
\proclaim{Lemma 3.1([Ma2;Theorem A and B], cf. [Ma5])}
Let $\frak L =(V, E, \lambda, \iota)$ be a $\lambda$-graph system over 
$\Sigma.$ 
Suppose that a $\lambda$-graph system $\frak L$ satisfies 
$\lambda$-condition (I).
Then the $C^*$-algebra ${\Cal O}_{\frak L}$ 
is the unique   
$C^*$-algebra generated by nonzero partial isometries
$s_\gamma, \gamma \in \Sigma$
and nonzero projections $e_i^l, i=1,2,\dots,m(l),\ l \in \Zp$  
satisfying  the following  operator relations:
$$
\align
& \sum_{\gamma \in \Sigma}  s_{\gamma}s_{\gamma}^*  =  1,\tag 3.1 \\
  \sum_{j=1}^{m(l)} e_j^l  & =  1, \qquad 
 e_i^l   =  \sum_{j=1}^{m(l+1)}I_{l,l+1}(i,j)e_j^{l+1},\tag 3.2 \\
& s_\gamma s_\gamma^* e_i^l =e_i^l 
s_\gamma s_\gamma^*, \tag 3.3 \\ 
s_{\gamma}^*e_i^l s_{\gamma} & =  
\sum_{j=1}^{m(l+1)} A_{l,l+1}(i,\gamma,j)e_j^{l+1},
\tag 3.4 
\endalign
$$
for
$
i=1,2,\dots,m(l),\l\in \Zp, 
 \gamma \in \Sigma,
 $
where $V_l = \{v_1^l,\dots, v_{m(l)}^l\}$ and 
$$
\align
A_{l,l+1}(i,\gamma,j)
 & =
\cases
1 &  
    \text{ if } \ s(e) = {v}_i^l, \lambda(e) = \gamma,
                       t(e) = {v}_j^{l+1} 
    \text{ for some }    e \in E_{l,l+1}, \\
0           & \text{ otherwise,}
\endcases \\
I_{l,l+1}(i,j)
 & =
\cases
1 &  
    \text{ if } \ \iota_{l,l+1}({v}_j^{l+1}) = {v}_i^l, \\
0           & \text{ otherwise}
\endcases 
\endalign
$$
for
$
i=1,2,\dots,m(l),\ j=1,2,\dots,m(l+1), \ \gamma \in \Sigma.
$ 
If in particular $\frak L$ is $\lambda$-irreducible,
the $C^*$-algebra ${\Cal O}_{\frak L}$
is simple and purely infinite.
\endproclaim
 We first consider the $C^*$-algebra ${\Cal O}_{\LCHLA}$ 
 for the $\lambda$-graph system $\LCHLA$.
\proclaim{Proposition 3.2}
Suppose that $A$ satisfies condition (I).
The $C^*$-algebra 
$
{\Cal O}_{\LCHLA} 
$ is canonically isomorphic to the Cuntz-Krieger algebra ${\Cal O}_A$ for the matrix $A$.
\endproclaim
\demo{Proof}
We notice that both the algebras 
$
{\Cal O}_{\LCHLA}
$
and
$
{\Cal O}_{A}
$
are uniquely determined by certain opertor relations
of their canonical generators.
We write the canonical generating partial isometries
and the projections in $
{\Cal O}_{\LCHLA}
$
as
$s_{\beta_i}, i=1,\dots,N$ and 
$e^l_{N( i_1 \cdots i_l)}, 
i_1\cdots i_l \in B_l(\Lambda_A), l \in \Zp$
respectively.
By the relations (3.1), (3.3) and (3.4), one has
$$
e^l_{N(i_1 \cdots i_l)}
=
\sum_{i_{l+1}, i_{l+2} =1}^N
s_{\beta_{i_1}}
e^{l+1}_{N(i_2 \cdots i_{l+1}i_{l+2})}
 s_{\beta_{i_1}}^*.
$$
For $l=1$, one sees that by (3.2)
$$
e_{N(i_1)}^1 
=
\sum_{i_2, i_3 =1}^N
s_{\beta_{i_1}}
e^2_{N(i_2 i_3)}
 s_{\beta_{i_1}}^*
 =
s_{\beta_{i_1}}
 s_{\beta_{i_1}}^*.
$$
As
$
\iota^2(v^{l+1}_{N(i_2\cdots i_{l+1}i_{l+2})})
=v^{l-1}_{N(i_2\cdots i_l)},
$
(3.2) implies the equality
$$
\sum_{i_{l+1}, i_{l+2} =1}^N
e^{l+1}_{N(i_2 \cdots i_{l+1}i_{l+2})}
= 
e^{l-1}_{N(i_2 \cdots i_l)}
$$
so that by induction one obtains 
$$
e^l_{N(i_1 \cdots i_l)}
=
s_{\beta_{i_1}}\cdots s_{\beta_{i_l}}
s_{\beta_{i_l}}^*\cdots s_{\beta_{i_1}}^*.
$$
One also sees that (3.4) implies the equality
$$
s_{\beta_i}^*s_{\beta_i} 
= \sum_{j=1}^N A(i, j) 
s_{\beta_j}s_{\beta_j}^*.
$$
As the equality
$
\sum_{i=1}^{N} s_{\beta_i}s_{\beta_i}^* = 1
$ 
holds,
 the $C^*$-algebra generated by partial isometries
$s_{\beta_i}, i=1,\dots,N$ is canonically isomorphic to 
the Cuntz-Krieger algebra ${\Cal O}_A$.
\qed
\enddemo 
In what follows, an $N \times N$ matrix 
$A$ is assumed to be irreducible  with entries in $\{0,1\}$, 
and satisfy condition (I).
By Proposition 2.5, the algebra $\OLA$ is simple and purely infinite. 
We will describe concrete operator relations among 
the canonical generators of the algebra  
$\OLA$.
Let $A_{l,l+1}, I_{l,l+1}$ be the matrices  as in Lemma 3.1 for the $\lambda$-graph system
$\LCHDA$. 
We denote by $m(l)$
the number of the vertex set
$V_l = \{ v^l_1, \dots, v_{m(l)}^l\}$
of $\LCHDA$.
Let $s_{\gamma}, \gamma \in \Sigma$ 
and $e_i^l, i=1,\dots,m(l), l \in \Zp$ 
be the canonical generating partial isometries and projections of $\OLA$.  
They satisfy the relations (3.1), (3.2), (3.3) and  (3.4)
for $\LCHDA$.   
Define the operators 
$
S_i,,
T_i,
i=1,\dots,N
$ by setting 
$$
S_i := s_{\alpha_i},\quad
T_i := s_{\beta_i}\quad
\quad
\text{ for } i=1,\dots,N.
$$
\proclaim{Proposition 3.3}
The operators 
$
S_i, T_i, i=1,\dots,N
$
satisfy the relations (1.1), (1.2), (1.3) and (1.4),
and generate the $C^*$-algebra $\OLA$.
\endproclaim
\demo{Proof}
The equality (1.1) is nothing but (3.1).
To prove (1.2), 
by the equality (3.4) and the first equality  of (3.2),
one has for a fixed $l \in \Zp$,
$$
\sum_{j=1}^N S_j^* S_j 
= \sum_{j=1}^N \sum_{i=1}^{m(l)}\sum_{k=1}^{m(l+1)}
        A_{l,l+1}(i,\alpha_j,k)e_k^{l+1}.
$$
For  $k = 1,\dots,m(l+1),$ 
there exists a unique edge labeled $\Sigma^-$ in
$\LCHDA$ whose terminal is
$v_k^{l+1}$.
Hence we have
$
\sum_{j=1}^{N}\sum_{i=1}^{m(l)}A_{l,l+1}(i,\alpha_j,k) = 1
$
so that 
$$
\sum_{j=1}^N S_j^* S_j 
= \sum_{k=1}^{m(l+1)} e_k^{l+1} = 1.
$$
For (1.3),
one similarly has
$$
T_i^* T_i 
= \sum_{k=1}^{m(l)}s_{\beta_i}^*e_k^l s_{\beta_i} 
= \sum_{h=1}^{m(l+1)} \sum_{k=1}^{m(l)} A_{l,l+1}(k,\beta_i,h) e_h^{l+1}.
$$
On the other hand,
$$
\sum_{j=1}^{N}A(i,j) S_j^* S_j
 = \sum_{h=1}^{m(l+1)}(\sum_{k=1}^{m(l)}\sum_{j=1}^{N}
A(i,j)
A_{l,l+1}(k,\alpha_j,h))e_h^{l+1}. 
$$
Let $h$ be written as 
$N(h_1\cdots h_{l+1})$.
Then the condition 
 $A_{l,l+1}(k,\beta_i,h) = 1$
is equivalent to the condition that
$i h_1 \in B_2(\Lambda_A)$
and
$k = N(i h_1\cdots h_{l-1})$.
On the other hand,
the condition
$
\sum_{j=1}^{N}
A(i,j)
A_{l,l+1}(k,\alpha_j,h) =1
$
is equivalent to the condition that
$j= h_1, A(i,j) =1$ for some $j$
and
$k = N(h_2\cdots h_{l+1})$.
Hence one has
$$
\sum_{k=1}^{m(l)} A_{l,l+1}(k,\beta_i,h)
=
\sum_{k=1}^{m(l)}\sum_{j=1}^{N}
A(i,j)
A_{l,l+1}(k,\alpha_j,h).
$$
This implies the equality (1.3).
For (1.4),
we put
$$
E_{\mu_1\cdots\mu_k} 
= S_{\mu_1}^*\cdots  S_{\mu_k}^* S_{\mu_k}\cdots S_{\mu_1}.
$$
By using the first equality of (3.2), (3.3) and (3.4) recursively,
one knows that the above defined projections commute with
$
S_jS_j^*
$
and
$T_j T_j^*$
for $j=1,\dots,N$.
Hence by (1.1), it follows that
$$
E_{\mu_1\cdots\mu_k} 
= \sum_{j=1}^{N} S_jS_j^* E_{\mu_1\cdots\mu_k}  S_jS_j^* + 
\sum_{j=1}^{N}T_jT_j^* E_{\mu_1\cdots\mu_k}   T_jT_j^*.
$$
As
$S_j^*E_{\mu_1\cdots\mu_k} S_j = A(j,\mu_1)S_j^*E_{\mu_1\cdots\mu_k} S_j $
and
$
S_{\mu_1}T_j = s_{\alpha_{\mu_1}} s_{\beta_j} = 0
$
if
$\mu_1 \ne j$,
one has 
 $$
E_{\mu_1\cdots\mu_k} 
 = \sum_{j=1}^{N} A(j,\mu_1)S_jS_j^* E_{\mu_1\cdots\mu_k}  S_jS_j^* + 
T_{\mu_1}T_{\mu_1}^* E_{\mu_1\cdots\mu_k}   T_{\mu_1}T_{\mu_1}^*.
$$
Since 
$A_{0,1}(1,\alpha_{\mu_1},j)=1$
if and only if 
$j = \mu_1,$
it follows that  by (3.4), 
$$
E_{\mu_1} 
= s_{\alpha_{\mu_1}}^*s_{\alpha_{\mu_1}} 
= \sum_{j=1}^{m(1)}A_{0,1}(1,\alpha_{\mu_1},j)e_j^{1}
= e_{\mu_1}^1.
$$
By (3.2) and (3.4), we similarly have
$$
\align
E_{\mu_1 \cdots \mu_k} 
& = s_{\alpha_{\mu_1}}^* \cdots s_{\alpha_{\mu_k}}^*
    s_{\alpha_{\mu_k}}\cdots s_{\alpha_{\mu_1}} \\
& =\sum_{{i_1}=1}^{m(1)}\cdots \sum_{{i_k}=1}^{m(k)}
A_{0,1}(1,\alpha_{\mu_k},i_1)\cdots A_{k-1,k}(i_{k-1},\alpha_{\mu_1},i_k)
e_{i_k}^{k}.
\endalign
$$
As 
$
\sum_{{i_1}=1}^{m(1)}\cdots \sum_{{i_{k-1}}=1}^{m(k-1)}
A_{0,1}(1,\alpha_{\mu_k},i_1) \cdots A_{k-1,k}(i_{k-1},\alpha_{\mu_1},i_k)
=1
$ 
if and only if
$i_k = N(\mu_1\cdots \mu_k),$
one knows 
$E_{\mu_1\cdots  \mu_k}= e_{N(\mu_1\cdots  \mu_k)}^k$.
Hence we have
$$
T_{\mu_1}^* E_{\mu_1\cdots  \mu_k} T_{\mu_1}
=
\sum_{j=1}^{m(k+1)} 
A_{k,k+1}(N(\mu_1\cdots \mu_k), \beta_{\mu_1},j) e_j^{k+1}.
$$
Since 
$A_{k,k+1}(N(\mu_1 \cdots \mu_k), \beta_{\mu_1},j)=1$
if and only if 
$
j= N(\mu_2 \cdots \mu_k \mu_{k+1} \mu_{k+2})
$
for some
$
\mu_{k+1}, \mu_{k+2} =1,\dots,N,
$
and the equality 
$$
\sum_{\mu_{k+1},\mu_{k+2}=1,\dots,N}
E_{{\mu_2}\cdots {\mu_k}{\mu_{k+1}}{\mu_{k+2}}}
=E_{{\mu_2}\cdots {\mu_k}}
$$
holds,
we have
$
T_{\mu_1}^* E_{\mu_1\cdots \mu_k}T_{\mu_1}= E_{\mu_2\cdots \mu_k}
$
so that 
$$
T_{\mu_1}T_{\mu_1}^* E_{\mu_1\cdots \mu_k}T_{\mu_1}T_{\mu_1}^*
= T_{\mu_1} E_{\mu_2\cdots \mu_k}T_{\mu_1}^*.
$$
Thus we conclude that (1.4) holds.
Consequently the operators 
$S_i, T_i,i=1,\dots,N$ 
satisfy the relations (1.1), (1.2), (1.3) and (1.4).   

In the above discussions, we have proved the equality
$$
e_{N(\mu_1\cdots  \mu_k)}^k= E_{\mu_1\cdots  \mu_k}
(= S_{\mu_1}^*\cdots S_{\mu_k}^*S_{\mu_k}\cdots  S_{\mu_1})
$$
for $\mu_1 \cdots \mu_k \in B_k(\Lambda_A)$.
Hence the algebra $\OLA$
is generated by
$S_1,\dots,S_N, $
$T_1,\dots,T_N$.
\qed
\enddemo

We next show that the relations (1.1), (1.2), (1.3) and (1.4) 
imply the relations (3.1), (3.2), (3.3) and (3.4).
Let $S_i, T_i, i=1,\dots,N$ 
be partial isometries  satisfying the operator relations
(1.1), (1.2), (1.3) and (1.4).
In the relation (1.4) for $k=2$, by summing up $\mu_2$ over
$\{1,\dots,N \}$ and using (1.2), 
we have      
$$
S_i^*S_i  =  \sum_{j=1}^N A(j,i) S_jS_j^*S_i^*S_iS_jS_j^* + T_iT_i^*,
\qquad i=1,\dots, N. 
\tag 3.5 
$$
\proclaim{Lemma 3.4} 
\roster
\item"(i)" 
$T_{i}^* S_j^* S_j T_{i}  
= 
\cases
T_i^* T_i & \text{ if } i=j,\\
0 &  \text{ if } i\ne j.
\endcases
$
\item"(ii)"
$
T_i^* E_{\mu_1 \cdots \mu_l} T_i = 
\cases
A(i,\mu_2)E_{\mu_2 \cdots \mu_l} & \text{ if } i=\mu_1,\\
0                                 & \text{ if } i\ne \mu_1\\
\endcases
$
\endroster
for $ l >1$, where
$E_{\mu_1 \cdots \mu_l} = S_{\mu_1}^*\cdots S_{\mu_l}^*S_{\mu_l}\cdots  S_{\mu_1}$
for $\mu_1 \cdots \mu_l \in B_l(\Lambda_A)$.
\endproclaim
\demo{Proof}
(i) By (3.5), we have
$$
T_i^*S_i^*S_i T_i = \sum_{j=1}^N A(j,i) T_i^*S_jS_j^*S_i^*S_iS_jS_j^* T_i
+ T_i^*T_iT_i^* T_i.
$$
The equality (1.1) implies $T_i^* S_j = 0$ for $i,j = 1,\dots,N$
and hence we have
$T_{i}^* S_i^* S_i T_{i}  = T_i^* T_i.$
 By (1.2), one has 
$$
\sum_{j=1}^{N}T_i^*S_j^* S_j T_i = T_i^* T_i
$$
so that
$
T_i^*S_j^* S_j T_i =0
$
for $i \ne j$.

(ii) By (1.4), we have
$$
T_i^* E_{\mu_1 \cdots \mu_l} T_i = 
\sum_{j=1}^N A(j,\mu_1)T_i^* S_jS_j^*
E_{\mu_1\cdots \mu_l}
S_jS_j^* T_i + T_i^*T_{\mu_1}E_{\mu_2\cdots \mu_l}T_{\mu_1}^*T_i
$$
for $l > 1$.
Since 
$ 
T_i^* S_j = 0
$
for
$ i,j= 1,\dots,N$
and
$
T_i^* T_{\mu_1} = 0 
$
for 
$ 
i\ne \mu_1,
$
we have
$$
T_i^* E_{\mu_1 \cdots \mu_l} T_i 
 = T_i^*T_{\mu_1}E_{\mu_2\cdots \mu_l}T_{\mu_1}^*T_i 
 =
{\cases 
T_i^*T_i E_{\mu_2\cdots \mu_l}T_i^*T_i & \text{ if } i = \mu_1, \\
0 & \text{ otherwise.}
\endcases} 
$$
By (1.3) one has 
$$
T_i^*T_i E_{\mu_2\cdots \mu_l}T_i^*T_i
=
\sum_{j=1}^{N}\sum_{k=1}^{N} 
A(i,j)A(i,k) S_j^*S_j S_{\mu_2}^* \cdots S_{\mu_l}^*S_{\mu_l}\cdots S_{\mu_2}
S_k^* S_k.
$$
By (1.2), one sees that
$
 S_j^*S_j S_{\mu_2}^* = 0 
$
for
$ 
j \ne \mu_2,
$ 
and 
$
 S_{\mu_2}^*S_{\mu_2} S_k^*S_k = 0 
$
for
$ k \ne \mu_2.
$
It then follows that
$$
T_i^* E_{\mu_1\cdots \mu_l} T_i
= A(i,\mu_2) E_{\mu_2\cdots \mu_l}.
$$
\qed
\enddemo
\proclaim{Lemma 3.5} 
Keep the above notations.
The projection
$E_{\mu_1 \cdots \mu_l}$ commutes with both $S_j S_j^*$ and $T_jT_j^*$
\endproclaim
\demo{Proof}
By (1.4), we have for $l > 1$
$$
 S_i S_i^* E_{\mu_1\cdots \mu_l} 
 = \sum_{j=1}^N A(j,\mu_1) S_i S_i^* S_jS_j^*
E_{\mu_1\cdots \mu_l}
S_jS_j^* + S_i S_i^* T_{\mu_1}E_{\mu_2\cdots \mu_l}T_{\mu_1}^*.
$$
By (1.1), one has
$
S_i^* T_{\mu_1} =0
$
for all
$ i,\mu_1,$
and by (1.2) one has 
$ 
S_i^*S_j =0 
$ for
$ i\ne j.
$
Hence we have
$
 S_i S_i^* E_{\mu_1\cdots \mu_l} 
= A(i,\mu_1) 
 S_i S_i^* E_{\mu_1\cdots \mu_l} 
S_iS_i^*
$ 
and similarly 
$
  E_{\mu_1\cdots \mu_l} S_i S_i^*
= A(i,\mu_1)  S_i S_i^* E_{\mu_1\cdots \mu_l} 
S_iS_i^*
$
so that 
$
 S_i S_i^*
$
commutes with
$ E_{\mu_1\cdots \mu_l}.$
By (1.4) and (1.1), we have
$$
\align
 T_i T_i^* E_{\mu_1\cdots \mu_l} 
& = \sum_{j=1}^N A(j,\mu_1) T_i T_i^* S_jS_j^*
E_{\mu_1\cdots \mu_l}
S_jS_j^* + T_i T_i^* T_{\mu_1}E_{\mu_2\cdots \mu_l}T_{\mu_1}^* \\
& = 
{\cases
T_{\mu_1} E_{\mu_2\cdots \mu_l} T_{\mu_1}^* & \text{ if } i =\mu_1,\\ 
 0  & \text{ otherwise. } \\
 \endcases}
\endalign
$$
We similarly have the same equality for 
$
  E_{\mu_1\cdots \mu_l} T_i T_i^*
$ 
as above
so that 
$ 
 T_i T_i^* 
$
commutes with
$E_{\mu_1\cdots \mu_l}.
$

For $l=1$, as we have 
$E_{\mu_1} = \sum_{\mu_2=1}^{N} E_{\mu_1 \mu_2}$ by (1.2),
the projection
$E_{\mu_1}$ commutes with both $S_j S_j^*$ and $T_j T_j^*$ by the above discussions.
\qed
\enddemo
\proclaim{Lemma 3.6}
Keep the above notations.
For $\mu_1, \dots,\mu_l \in \{1,\dots,N\}$ we have
$E_{\mu_1\cdots\mu_l} = 0$ if $\mu_1\cdots\mu_l \not\in B_l(\Lambda_A).$ 
\endproclaim
\demo{Proof}
As we are assuming that the matrix $A$ has no zero rows or columns,
for $l=1$, $B_1(\Lambda_A) = \{ 1,\dots,N\}$.
By (3.5) one has  for $\mu_1 = 1,\dots,N$
$$ 
S_{\mu_1}^*S_{\mu_1}  =  \sum_{j=1}^N A(j,\mu_1) S_jS_j^*S_{\mu_1}^*S_{\mu_1}
S_jS_j^* + T_{\mu_1}T_{\mu_1}^*
$$
so that for $\mu_0 = 1,\dots,N$
$$
S_{\mu_0}^*S_{\mu_1}^*S_{\mu_1}S_{\mu_0}  
=   A(\mu_0,\mu_1) S_{\mu_0}^*S_{\mu_1}^*S_{\mu_1}S_{\mu_0}
$$
because $S_{\mu_0}^*S_j =0$ if $\mu_0 \ne j$, and
$S_{\mu_0}^*T_{\mu_1} = 0$ by (1.1).
This means that 
$E_{\mu_0\mu_1} = 0$ if $\mu_0\mu_1 \not\in B_2(\Lambda_A).$ 

Suppose next that the assertion holds for $l=k >1$.
By (1.4) one has for $\mu_1\cdots\mu_k \in B_k(\Lambda_A)$ and $\mu_0 = 1,\dots,N$
$$
S_{\mu_0}^*
E_{\mu_1\cdots \mu_k} S_{\mu_0}= \sum_{j=1}^N A(j,\mu_1)S_{\mu_0}^* S_jS_j^*
E_{\mu_1\cdots \mu_k}
S_jS_j^* S_{\mu_0}
+ S_{\mu_0}^*T_{\mu_1}E_{\mu_2\cdots \mu_k}T_{\mu_1}^*S_{\mu_0}
$$
so that we have
$$
E_{\mu_0 \mu_1 \cdots \mu_k} = A(\mu_0,\mu_1)E_{\mu_0 \mu_1 \cdots \mu_k}.
$$
For
$\mu_1 \cdots \mu_k \in B_k(\Lambda_A)$,
we have
$\mu_0 \mu_1 \cdots \mu_k \not\in B_{k+1}(\Lambda_A) $
if and only if $A(\mu_0,\mu_1) =0$.
Hence the assertion holds for $l=k+1$.
Therefore the assertion holds for all $l \in \Zp$.
\qed
\enddemo
\proclaim{Proposition 3.7} 
Keep the above notations.
Put
$$
 s_{\alpha_i}:=S_i ,\quad
  s_{\beta_i}:=T_i \quad \text{ for }
 \quad i=1, \dots, N 
$$
and
$$
\align
e^0_1 & := 1,\\
e_{N(\mu_1 \cdots \mu_l)}^l & := E_{\mu_1 \cdots \mu_l}( =
S_{\mu_1}^*\cdots S_{\mu_l}^* S_{\mu_l}\cdots S_{\mu_1})
\qquad \text{ for } 
\mu = \mu_1\cdots \mu_l\in B_l(\Lambda_A).\\
\endalign
$$
Then the family 
of operators
$s_\gamma, \gamma \in \Sigma,$ 
$e_{N(\mu_1 \cdots \mu_l)}^l,\mu_1\cdots \mu_l\in B_l(\Lambda_A)$  
satisfies the operator relations
(3.1), (3.2), (3.3) and (3.4) for the $\lambda$-graph system $\LCHDA$.
\endproclaim
\demo{Proof}
The relation (3.1) is nothing but the equality (1.1).
The equality (1.2) implies 
$
\sum_{\mu_1\in B_1(\Lambda_A)}e_{N(\mu_1)}^1 = 1.
$ 
Suppose that
$\sum_{\mu_1\cdots \mu_l\in B_l(\Lambda_A)}e_{N(\mu_1 \cdots \mu_l)}^l =1$ 
holds for $l=k$.
As
$$
S_{\mu_1}^*\cdots S_{\mu_k}^*S_{\mu_k}\cdots S_{\mu_1} 
=  \sum_{h=1}^N 
S_{\mu_1}^*\cdots S_{\mu_k}^*
S_h^*S_h 
S_{\mu_k}\cdots S_{\mu_1},
$$
the equality
$\sum_{\mu_1\cdots \mu_l\in B_l(\Lambda_A)}e_{N(\mu_1 \cdots \mu_l)}^l =1$
 holds for $l= k+1$ by Lemma 3.6 and hence for all $l$.
The above equality with the equality
$$
I_{l,l+1}(N(\mu_1\cdots\mu_l), N(\nu_1\cdots\nu_{l+1}))
=
\cases
1 & \text{ if } \nu_1 \cdots \nu_l = \mu_1 \cdots \mu_l,\\
0 & \text{ otherwise } 
\endcases
$$
for $\nu_1\cdots\nu_{l+1} \in B_{l+1}(\Lambda_A)$
implies the second relation of (3.2) by using Lemma 3.6.
The equality
(3.3) comes from  Lemma 3.5.

We will finally show the equality  (3.4).

For $l=0$, $e^0_1 = 1$ by definition.
If $\gamma =\alpha_k$ for some $k=1,\dots,N$,
one has $A_{0,1}(1,\alpha_k,j) =1$
if and only if $j=k$.
Hence
$$
s_{\alpha_k}^*e_1^0 s_{\alpha_k}  =  s_{\alpha_k}^* s_{\alpha_k}
= e^1_k =
\sum_{j=1}^{m(1)} A_{0,1}(1,\alpha_k,j)e_j^{1}.
$$
If $\gamma =\beta_k$ for some $k=1,\dots,N$,
one has $A_{0,1}(1,\beta_k,j) =A(k,j).$
Hence by the relation (1.3)
one has
$$
s_{\beta_k}^*e_1^0 s_{\beta_k}  =  T_k^* T_k
= e^1_k =
\sum_{j=1}^{N} A_{0,1}(k,j)S_j^* S_j
=
\sum_{j=1}^{N} A_{0,1}(1,\beta_k,j)e_j^{1}.
$$

For $l=1$, one sees that
$e^1_i = e^1_{N(i)}$.
If $\gamma =\alpha_k$ for some $k=1,\dots,N$,
one has 
$$
s_{\alpha_k}^*e_i^1 s_{\alpha_k}  =  S_k^* S_i^* S_i S_k
= e^2_{N(ki)} =
\sum_{j=1}^{m(2)} A_{1,2}(i,\alpha_k,j)e_j^{2},$$
where the last equality above comes from the fact that 
$A_{1,2}(i,\alpha_k,j) =1$
if and only if $j=N(ki)$.
If $\gamma =\beta_k$ for some $k=1,\dots,N$,
one has by Lemma 3.4 (i) and (1.3)
$$
s_{\beta_k}^*e_i^1 s_{\beta_k}  =  T_k^*S_i^* S_i T_k
= 
\cases
\sum_{j=1}^N A(i,j)S_j^*S_j & \text{ if } k=i,\\
0 & \text{ if } k\ne i.
\endcases
$$
By Lemma 3.6 and (1.2), one has
$$
\sum_{j=1}^N A(i,j) S_j^* S_j 
= \sum_{\mu_1 \mu_2 \in B_2(\Lambda_A)}A(i,\mu_1)A(\mu_1,\mu_2) 
S_{\mu_1}^*S_{\mu_2}^*S_{\mu_2}S_{\mu_1}.
$$
Since
$
A_{1,2}(i,\beta_k,\mu_1\mu_2) =1
$
if and only if 
$k=i,$
$A(i,\mu_1)=A(\mu_1,\mu_2)=1$,
it follows that  by
$
S_{\mu_1}^*S_{\mu_2}^*S_{\mu_2}S_{\mu_1}
=
e_{N(\mu_1\mu_2)}^2,
$
$$
s_{\beta_k}^*e_i^1 s_{\beta_k}  =  
\sum_{j=1}^{m(2)}A_{1,2}(i,\beta_k,j) e_j^2.
$$

For $\mu_1\cdots\mu_l \in B_l(\Lambda_A)$ with $l >1$
and 
$\alpha_k \in \Sigma^-,$
the relation (1.4)
implies
$$
\align
s_{\alpha_k}^* e_{N(\mu_1\cdots\mu_l) }^l s_{\alpha_k} & 
= A(k,\mu_1)
S_{k}^*
S_{\mu_{1}}^*\cdots S_{\mu_l}^*S_{\mu_l}\cdots S_{\mu_{1}}
S_{k}\\
& = A_{l,l+1}(N(\mu_1\cdots\mu_l) ,\alpha_k, N(k \mu_1\cdots\mu_l) )
e_{N(k \mu_1\cdots\mu_l) }^{l+1}.
\endalign
$$
Since 
$A_{l,l+1}(N(\mu_1\cdots \mu_l), \alpha_k,i) = 0$ 
if $i \ne N(k\mu_1\cdots \mu_l)$,
one has
$$
s_{\alpha_k}^* e_{N(\mu_1\cdots \mu_l)}^l s_{\alpha_k}
= \sum_{\nu_1 \cdots \nu_{l+1}\in B_{l+1}(\Lambda_A)}
A_{l,l+1}(N(\mu_1\cdots \mu_l),\alpha_k, N(\nu_1 \cdots \nu_{l+1})) 
e_{N(\nu_1 \cdots \nu_{l+1})}^{l+1}.
$$
We also have 
by Lemma 3.4
$$
\align
& s_{\beta_j}^* e_{N(\mu_1\cdots \mu_l)}^l s_{\beta_j} \\
& = 
T_{j}^*
E_{\mu_1 \cdots \mu_l} T_{j}\\
& 
=
\cases
A(j,\mu_2) E_{\mu_2 \cdots \mu_l}
 & \text{ if } j= \mu_1,\\
0 & \text{ otherwise} 
\endcases \\
& 
=
\cases
\sum \Sb \mu_{l+1}, \mu_{l+2} \in B_2(\Lambda_A) \endSb
A(j,\mu_2) S_{\mu_{2}}^*\cdots S_{\mu_l}^*  
S_{\mu_{l+1}}^*S_{\mu_{l+2}}^*
S_{\mu_{l+2}}S_{\mu_{l+1}}
S_{\mu_l}\cdots S_{\mu_{2}}
 & \text{ if } j= \mu_1,\\
0 & \text{ otherwise.} 
\endcases 
\endalign
$$
Since 
one has
$$
\align
& A_{l,l+1}(N(\mu_1\cdots \mu_l), \beta_j, N(\nu_1\cdots \nu_{l+1})) \\
=
& {\cases
1 & \text{ if } j=\mu_1,\,  A(j,\mu_2) = 1
\text{ and } \nu_i = \mu_{i+1} \text{ for } i=1,\dots,l-1,\\
0 & \text{ otherwise,}
\endcases}
\endalign
$$
we have
$$
s_{\beta_j}^* e_{N(\mu_1\cdots \mu_l)}^l s_{\beta_j}
=\sum_{\nu_1\cdots \nu_{l+1} \in B_{l+1}(\Lambda_A)}  
A_{l,l+1}(N(\mu_1\cdots \mu_l),\beta_j, N(\nu_1\cdots\nu_{l+1}))
e_{N(\nu_1\cdots\nu_{l+1})}^{l+1}.
$$
Therefore (3.4) holds
\qed
\enddemo
By a general theory of the $C^*$-algebras associated with $\lambda$-graph systems [Ma2], the algebras $\OLA$ are nuclear.
If $A$ is irreducible with condition (I) in the sense of Cuntz-Krieger,
the $\lambda$-graph system $\LCHDA$ is $\lambda$-irreducible with condition (I) from Proposition 2.5 so that the algebra
$\OLA$ is simple and purely infinite by [Ma5].  
By Proposition 3.3 and Proposition 3.7, 
the family of the operator relations (1.1), (1.2), 
(1.3) and  (1.4) is equivalent to 
the family of the operator relations  (3.1), 
(3.2), (3.3) and  (3.4).
Thus by Lemma 3.1 
we conclude Theorem 1.1.

\heading 4. K-Theory
\endheading
In this section, we will present K-theory formulae of the $C^*$-algebra 
$\OLA$ in terms of the topological Markov shift defined by the matrix $A$.
Recall that the right one-sided topological Markov shift 
$(X_{\Lambda_A}, \sigma_A)$
is a continuous map $\sigma_A$ on $X_{\Lambda_A}$ where 
$$
\align
X_{\Lambda_A}& = \{ {(x_i)}_{i \in \Bbb N}\in \{ 1,\dots,N \}^{\Bbb N}
 \mid A(x_i,x_{i+1}) = 1, i \in \Bbb N \},\\
\sigma_A({(x_i)}_{i \in \Bbb N})& = {(x_{i+1})}_{i \in \Bbb N}
\qquad \text{ for } {(x_i)}_{i \in \Bbb N} \in X_{\Lambda_A}.
\endalign
$$ 
The space $X_{\Lambda_A}$ is naturally identified with 
$X_{D_A^+}$ in the proof of Proposition 2.1. 
Let
$C(X_{\Lambda_A},\Bbb Z)$ be the abelian group of all $\Bbb Z$-valued continuous functions on $X_{\Lambda_A}$.
Define endomorphisms $\sigma_{\Lambda_A}$ and $\lambda_{\Lambda_A}$ on
$C(X_{\Lambda_A},\Bbb Z)$ by
$$
\sigma_{\Lambda_A}(f) (x) = f(\sigma_{\Lambda_A}(x)),
\qquad
\lambda_{\Lambda_A}(f) (x) = \sum_{j=1}^Nf(jx) 
\qquad
\text{ for } 
f \in  C(X_{\Lambda_A},\Bbb Z)
$$
and
$ x= {(x_i)}_{i \in \Bbb N} \in X_{\Lambda_A},
$
where
$jx = (j,x_1,x_2, \dots) \in X_{\Lambda_A}$.
Let $S_i, T_i, i=1,\dots,N$ be the generating partial isometries 
of the $C^*$-algebra 
$\OLA$ as in Theorem 1.1.
Let $\ALA$ be the $C^*$-subalgebra of $\OLA$ generated by the projections
$
E_{\mu_1\cdots\mu_l}
= S_{\mu_1}^*\cdots S_{\mu_l}^*S_{\mu_l}\cdots S_{\mu_1},
\,
\mu_1\cdots\mu_l \in \Lambda_A^*$.
Define two endomorphisms $\lambda_{\Sigma^-}$ and  $\lambda_{\Sigma^+}$
on it by 
$$
\lambda_{\Sigma^-}(a)  = \sum_{j=1}^N S_j^* a S_j, \qquad   
\lambda_{\Sigma^+}(a)  = \sum_{j=1}^N T_j^* a T_j \qquad 
\text{ for } a \in \ALA
$$
Let $C(X_{\Lambda_A},\Bbb C)$ be the abelian $C^*$-algebra of all $\Bbb C$-valued continuous functions on $X_{\Lambda_A}$.
We note that its $K_0$-group $K_0(C(X_{\Lambda_A},\Bbb C))$
is naturally identified with $C(X_{\Lambda_A},\Bbb Z)$.
\proclaim{Lemma 4.1}
Let
$ \varPhi: \ALA \longrightarrow C(X_{\Lambda_A},\Bbb C)$ 
be a map defined by 
$$
\varPhi (E_{\mu_1\cdots\mu_l}) = \chi_{\mu_1\cdots\mu_l} \qquad \text{ for }
\mu_1\cdots\mu_l \in \Lambda_A^*
$$
where 
$
\chi_{\mu_1\cdots\mu_l}
$
is the characteristic function for the word 
$\mu_1\cdots\mu_l $ on $X_{\Lambda_A}$
defined by
$$
\chi_{\mu_1\cdots\mu_l} ({(x_i)}_{i \in \Bbb N}) = 
\cases
1 & \text{ if }  (x_1, \dots, x_l) = (\mu_1,\dots,\mu_l),\\
0 & \text{ otherwise}.
\endcases
$$
Then we have
\roster
\item"(i)" $\varPhi$ gives rise to an isomorphism from 
$ \ALA $ onto $C(X_{\Lambda_A},\Bbb C)$.
\item"(ii)"
Both of the diagrams
 $$
\CD
K_0(\ALA)
@>
\varPhi_* 
>> 
C(X_{\Lambda_A},\Bbb Z)
 \\
@V
\lambda_{\Sigma^-_*}  
VV @V
\sigma_{\Lambda_A} 
VV \\
K_0(\ALA)
@>\varPhi_* >> 
C(X_{\Lambda_A},\Bbb Z)
\endCD,
\qquad
\CD
K_0(\ALA)
@>
\varPhi_* 
>> 
C(X_{\Lambda_A},\Bbb Z)
 \\
@V
\lambda_{\Sigma^+_*}  
VV @V
\lambda_{\Lambda_A} 
VV \\
K_0(\ALA)
@>\varPhi_* >> 
C(X_{\Lambda_A},\Bbb Z)
\endCD
$$
are commutative,
where $\varPhi_*$ is the induced isomorphism from
$K_0(\ALA)$ to $K_0(C(X_{\Lambda_A},\Bbb C)) (=C(X_{\Lambda_A},\Bbb Z))$,
and $\lambda_{\Sigma^-_*}$, $\lambda_{\Sigma^+_*}$ are induced endomorphisms
on $K_0(\ALA)$ by 
$\lambda_{\Sigma^-}$, $\lambda_{\Sigma^+}$ respectively.
\endroster
\endproclaim
\demo{Proof}
(i) The assertion is straightforward.

(ii)
The equality 
$$
\varPhi(\lambda_{\Sigma^-} (E_{\mu_1\cdots\mu_l}))
= \sum_{j=1}^N \chi_{j \mu_1\cdots\mu_l}
$$
is immediate. As
$$
\sigma_{\Lambda_A}(\chi_{\mu_1\cdots\mu_l}) (x) 
 =
{\cases
1 & \text{ if } (x_2,\dots,x_{l+1}) = (\mu_1,\dots,\mu_l),\\
0 & \text{ otherwise,}
\endcases}
$$
for 
$x = {(x_i)}_{i \in \Bbb N} \in X_{\Lambda_A}$,
the equality 
$$
\sigma_{\Lambda_A}(\chi_{\mu_1\cdots\mu_l}) = 
\varPhi(\lambda_{\Sigma^-} (E_{\mu_1\cdots\mu_l}))
$$
is clear.
Hence the first diagram is commutative.

For the second diagram, as 
$$
T_j^* E_{\mu_1\cdots\mu_l} T_j
 = 
\cases
A(j,\mu_2)E_{\mu_2 \cdots \mu_l} & \text{ if } j=\mu_1,\\
0                                 & \text{ if } j\ne \mu_1\\
\endcases
$$
by Lemma 3.4,
it follows  that 
$$
\align
\varPhi( T_j^* E_{\mu_1\cdots\mu_l} T_j)(x)
& = 
{\cases
A(j,\mu_2)\chi_{\mu_2 \cdots \mu_l}(x) & \text{ if } j=\mu_1,\\
0                                 & \text{ if } j\ne \mu_1\\
\endcases}\\
& = 
{\cases
1 & \text{ if } (\mu_1,\dots, \mu_l) =(j,x_1, x_2, \dots,x_{l-1}), \\
0 & \text{ otherwise.}
\endcases}\\
\endalign
$$
On the other hand, one sees for 
$ x = {(x_i)}_{i \in \Bbb N} \in X_{\Lambda_A}$
$$
\align
\lambda_{\Lambda_A}(\chi_{\mu_1\cdots\mu_l}) (x) 
& =
 \sum_{j=1}^N \chi_{\mu_1\cdots\mu_l} (j x) \\
& =
{\cases
1 & \text{ if }  (\mu_1,\dots,\mu_l) = (j, x_1,x_2,\dots,x_{l-1}) \text{ for some } j=1,\dots,N\\
0 & \text{ otherwise}
\endcases}
\endalign
$$
so that one obtains
$$
\lambda_{\Lambda_A}(\chi_{\mu_1\cdots\mu_l}) = 
\sum_{j=1}^N \varPhi(T_j^* E_{\mu_1\cdots\mu_l} T_j)
=\varPhi(\lambda_{\Sigma^+}(E_{\mu_1\cdots\mu_l})).
$$
Hence the second diagram is commtative.
\qed
\enddemo 
Therefore we have
\proclaim{Theorem 4.2}
\roster
\item"(i)"
$K_0(\OLA) = C(X_{\Lambda_A},\Bbb Z) 
/ (\id - (\sigma_{\Lambda_A} + \lambda_{\Lambda_A})) C(X_{\Lambda_A},\Bbb Z) .$
\item"(ii)"
$K_1(\OLA) = \Ker (\id - (\sigma_{\Lambda_A} + \lambda_{\Lambda_A}))
\quad \text{ in } C(X_{\Lambda_A},\Bbb Z) .$
\endroster
\endproclaim
\demo{Proof}
(i) By an argument of [Ma2;Theorem 5.5], one knows
$$
K_0(\OLA) = K_0(\ALA) 
/ (\id - {\lambda_{\LCHDA}}_*) K_0(\ALA).
$$
where
${\lambda_{\LCHDA}}_*$ is an endomorphism on 
$K_0(\ALA)$ induced by the map
$\lambda_{\LCHDA}: \ALA \rightarrow \ALA$
defined by
$$
\lambda_{\LCHDA}(a) = 
\sum_{\gamma \in \Sigma^- \cup \Sigma^+}S_\gamma^* a S_\gamma
\qquad \text{ for } a \in \ALA.
$$
As 
$\lambda_{\LCHDA}(a) = 
\lambda_{\Sigma^-}(a) + \lambda_{\Sigma^+}(a)
$,
one sees the desired formula by the previous lemma.

(ii)
Similarly 
by an argument of [Ma2;Theorem 5.5], one knows
$$
K_1(\OLA) = \Ker (\id - {\lambda_{\LCHDA}}_*) \text{ in }  K_0(\ALA)
$$
so that one sees the desired formula by the previous lemma.
\qed
\enddemo
We remark that the space $X_{\Lambda_A}$ is homeomorphic to a 
a Cantor discontinuum $\frak K$.

The formulae of Theorem 4.2 may be rewritten 
in terms of nonnegative matrix system
$(M_{l,l+1}, I_{l,l+1})_{l\in \Zp}$
for the $\lambda$-graph systems $\LCHDA$ (cf.[Ma]),
where
$I_{l,l+1}$ is 
the matrix as in Lemma 3.1 for $\LCHDA$ and
$M_{l,l+1}$ is 
the matrix defined by
$$
M_{l,l+1}(i,j) = \sum_{\gamma \in \Sigma}A_{l,l+1}(i,\gamma,j)
\quad
\text{ for }
i=1,\dots,m(l), \, j=1,\dots,m(l+1) 
$$
for the matrix $A_{l,l+1}$ in Lemma 3.1 for $\LCHDA$.
Then the relations 
$
I_{l,l+1} M_{l+1,l+2} = M_{l,l+1}I_{l+1,l+2}
$ 
for 
$
l \in \Zp
$
hold.
 The groups
$K_0(\OLA)$ and $ K_1(\OLA)$
are computed by the following formulae. 
\proclaim{Proposition 4.3([Ma2], cf. [Ma])}
\roster
\item"(i)"
$
K_0(\OLA ) 
=  
\underset{l}\to{\varinjlim} 
\{ {\Bbb Z}^{m(l+1)} / (M_{l,l+1}^t - I_{l,l+1}^t){\Bbb Z^{m(l)}}, 
   \bar{I}^t_{l,l+1} \},
$
where the inductive limit is taken along the natural induced homomorphisms
$
   \bar{I}^t_{l,l+1}, l \in \Zp
$
by the 
matrices $I^t_{l,l+1}$. 
\item"(ii)"
$
K_1(\OLA)
=\underset{l}\to{\varinjlim}\{ \Ker (M_{l,l+1}^t - I_{l,l+1}^t)
             \text{ in }  {\Bbb Z}^{m(l)}, I^t_{l,l+1} \}, 
$
where the inductive limit is taken along the homomorphisms of the restrictions of $I^t_{l,l+1}$ to  $\Ker (M_{l,l+1}^t - I_{l,l+1}^t).$
\endroster
\endproclaim

\heading 5. Examples
\endheading

\proclaim{Example 1 (Dyck shifts)}
\endproclaim
For the matrix $A$ all of whose entries are $1$, Theorem 1.1 goes to 
\proclaim{Proposition 5.1 ([Ma6])}
The $C^*$-algebra ${\Cal O}_{{\frak L}^{Ch(D_N)}}$
associated with the Cantor horizon $\lambda$-graph system
${\frak L}^{Ch(D_N)}$ for the Dyck shift $D_N$
is unital, separable,  nuclear, simple and purely infinite.
It is the unique $C^*$-algebra generated by $N$ partial isometries 
$S_i,  i=1,\dots,N$ and  $N$ isometries 
$ T_i, i=1,\dots,N$  
subject to the following operator relations:
$$
\sum_{j=1}^{N}   S_j^*S_j     = 1, \qquad
E_{\mu_1\cdots \mu_k} = \sum_{j=1}^N  S_jS_j^*
E_{\mu_1\cdots \mu_k}
S_jS_j^* + T_{\mu_1}E_{\mu_2\cdots \mu_k}T_{\mu_1}^*    
$$
where 
$
E_{\mu_1\cdots \mu_k}
= S_{\mu_1}^*\cdots S_{\mu_k}^*S_{\mu_k}\cdots S_{\mu_1},
$
$\mu_1,\dots,\mu_k\in \{ 1,\dots, N\}$.
The K-groups are 
$$
K_0({\Cal O}_{\frak L^{Ch(D_N)}}) \cong {\Bbb Z}/N{\Bbb Z} \oplus C(\frak K,\Bbb Z),
\qquad
K_1({\Cal O}_{\frak L^{Ch(D_N)}}) \cong 0.
$$
\endproclaim
\demo{Proof}
The relation (1.3) implies that $T_i, i=1,\dots, N$
are isometries.
By summing up $\mu_2$ over $\{ 1, \dots,N\}$
 in the second relation above for $k=2$,
 one has the equalities
$$
S_i^*S_i = \sum_{j=1}^N  S_jS_j^*S_i^*S_iS_jS_j^* + T_iT_i^*, 
\qquad i=1,\dots,N
$$
by using the first relation above.
By summing up $i=1,2,\dots,N$ in the above equalities,
one sees
the relation (1.1).
\qed
\enddemo

\proclaim{Example 2 (Fibonacci Dyck shift)}
\endproclaim

Let $F$ be the $2\times 2$ matrix
$
\left[\smallmatrix 
1 & 1 \\
1 & 0 \\
\endsmallmatrix
\right]
$.
It is  the smallest matrix in the irreducible square matices with condition (I) such that the associated topological Markov shift $\Lambda_F$
is not conjugate to any full shift.
The topological entropy of 
$\Lambda_F$  is
$\log \frac{1 + \sqrt{5}}{2}$ 
the logarithm of the Perron eigenvalue of $F$.     
We call the subshift $D_F$ the Fibonacci Dyck shift.
As the matrix is irreducible with condition (I),
the associated $C^*$-algebra
$\OLF$ is simple and purely infinite.

\proclaim{Proposition 5.2}
The $C^*$-algebra $\OLF$
associated with the $\lambda$-graph system
$\LCHDF$ is unital, separable, nuclear, simple and purely infinite.
It is the unique $C^*$-algebra generated by one isometry $T_1$
and three  partial isometries 
$S_1, S_2, T_2$  
subject to the following operator relations:
$$
\align
\sum_{j=1}^{2} & ( S_jS_j^* +   T_jT_j^* ) 
  =  \sum_{j=1}^{2}  S_j^*S_j      = 1,  \qquad 
T_2^*T_2   =   S_1^*S_1,  \\
E_{\mu_1\cdots \mu_k} & = \sum_{j=1}^2 F(j,\mu_1) S_jS_j^*
E_{\mu_1\cdots \mu_k}
S_jS_j^* + T_{\mu_1}E_{\mu_2\cdots \mu_k}T_{\mu_1}^*, \qquad k >1    
\endalign
$$
where 
$E_{\mu_1\cdots \mu_k}
= S_{\mu_1}^*\cdots S_{\mu_k}^*S_{\mu_k}\cdots S_{\mu_1}$,
$(\mu_1,\cdots,\mu_k)\in \Lambda_F^*,$ 
and $\Lambda_F^*$ is the set of admissible words of the topological Markov shift $\Lambda_F$ defined by the matrix $F$.
The K-groups are
$$
K_0(\OLF) \cong {\Bbb Z} \oplus  C({\frak K},\Bbb Z)^{\infty},
\qquad
K_1(\OLF) \cong 0.
$$
\endproclaim
\demo{Proof}
The operator relations above  directly come from Theorem 1.1.
The K-group formulae above are not direct.
Its computations need some tecnichal steps as in [Ma7].
The full proof of the above K-group formulae are written in [Ma7].
\qed
\enddemo


\bigskip
{\it e-mail}: kengo{\@}yokohama-cu.ac.jp

\bye